\documentclass[12pt]{amsart}
\usepackage{mathrsfs}
\usepackage{amsfonts}
\usepackage{amssymb}
\usepackage{amsfonts, amscd, amsmath, mathrsfs, amssymb, amsthm, amsxtra, bbding, epsfig, graphicx, latexsym, url, mathbbol, bbold}
\usepackage[papersize={7.7in,10.3in},textwidth=14.8cm,textheight=20.5cm,centering]{geometry}
\usepackage{enumerate}

\usepackage{xcolor}
\definecolor{cite}{rgb}{0.00,0.60,1.00}
\definecolor{url}{rgb}{1.00,0.10,0.80}
\definecolor{link}{rgb}{0.00,0.00,1.00}
\usepackage[colorlinks,linkcolor=link,urlcolor=url,citecolor=cite,pagebackref,breaklinks]{hyperref}

\hypersetup{
pdfstartpage=1,
pdfstartview=FitH}

\DeclareFontFamily{U}{mathx}{\hyphenchar\font45}
\DeclareFontShape{U}{mathx}{m}{n}{
      <5> <6> <7> <8> <9> <10>
      <10.95> <12> <14.4> <17.28> <20.74> <24.88>
      mathx10
      }{}
\DeclareSymbolFont{mathx}{U}{mathx}{m}{n}
\DeclareMathAccent{\widecheck}{\mathalpha}{mathx}{"71}

 \usepackage{caption} 
\numberwithin{equation}{section}

\allowdisplaybreaks

\newtheorem{theorem}{Theorem}[section]
\newtheorem{lemma}{Lemma}[section]

\newtheorem{corollary}{Corollary}[section]

\makeatletter
\newcounter{roem}
\renewcommand{\theroem}{\Roman{roem}}

\newcommand{\c@org@eq}{}
\let\c@org@eq\c@equation
\newcommand{\org@theeq}{}
\let\org@theeq\theequation

\newcommand{\setroem}{
\let\c@equation\c@roem
 \let\theequation\theroem}

\newcommand{\setarab}{
\let\c@equation\c@org@eq
\let\theequation\org@theeq}
\makeatother

\newtheorem*{claim*}{Claim}

\theoremstyle{remark}

\newcommand{\Tr}{\mathrm{Tr}}

\newcommand{\GL}{\mathrm{GL}}

\newcommand{\PGL}{\mathrm{PGL}}

\DeclareMathOperator{\Mod}{mod}

\renewcommand{\bmod}[1]{\,(\Mod{ #1})}

\newcommand{\bx}{\mathbf{x}}
\newcommand{\by}{\mathbf{y}}

\newcommand{\bF}{\mathbf{F}}

\newcommand{\cH}{\mathcal{H}}

\newcommand{\fn}{\mathfrak{n}}

\newcommand{\ft}{\mathfrak{t}}

\newcommand{\fG}{\mathfrak{G}}

\usepackage{graphicx}
\usepackage{tikz}

\begin{document}

\title{A double character sum of Conrey--Iwaniec and Petrow--Young}

\author{Ping Xi}

\address{School of Mathematics and Statistics, Xi'an Jiaotong University, Xi'an 710049, P. R. China}

\email{ping.xi@xjtu.edu.cn}

\begin{abstract} 
We show that a double character sum, appearing in the work of Conrey--Iwaniec and Petrow--Young on Weyl bound for certain $L$-functions, is essentially a hypergeometric sum introduced by Katz. This produces a simple proof of the upper bound for this sum.
\end{abstract}



\maketitle

\setcounter{tocdepth}{1}

\section{Introduction}\label{sec:Introduction}

In their deep studies on the subconvexity of automorphic $L$-functions, Conrey and Iwaniec \cite{CI00} introduced the following character sum
\begin{align*}
g(\chi,\eta)=\mathop{\sum\sum}_{u,v\bmod p}\chi(u)\overline{\chi}(u+1)\overline{\chi}(v)\chi(v+1)\eta(uv-1),
\end{align*}
where $p$ is an odd prime, $\chi,\eta$ are multiplicative characters modulo $p$. They proved that 
\begin{align}\label{eq:g-upperbound}
g(\chi,\eta)\ll p
\end{align}
as long as $\chi$ is quadratic and $\eta$ is non-trivial (in which case $g(\chi,\eta)$ is real-valued).
The bound \eqref{eq:g-upperbound} was recently generalized to all non-trivial characters $\chi,\eta$ by
Petrow and Young \cite{PY20} when they extend the subconvexity of Conrey and Iwaniec \cite{CI00} to characters of any order $>1$. Both proofs utilize Deligne's work on Riemann Hypothesis for varieties over finite fields. See also \cite{FKM21} for another elegant proof in the language of their cohomological transforms.

In this note, we show that $g(\chi,\eta)$ is essentially a hypergeometic sum introduced by Katz, so that the desired upper bound for the former one then follows from that of the latter.
To proceed, we would like to work in finite field settings. Let $k$ be a finite field with $q$ elements of characteristic $p$. For $\chi,\eta$ non-trivial characters of $k^\times,$ 
define the double character sum
\begin{align}\label{eq:charactersum}
g(\chi,\eta;k)=\mathop{\sum\sum}_{u,v\in k}\chi\Big(\frac{u(v+1)}{v(u+1)}\Big)\eta(uv-1).
\end{align}
We now quickly state our result.
\begin{theorem}\label{thm:g-hyper}
Let $\chi,\eta$ be two non-trivial characters in $k^\times.$ Then we have
\begin{align*}
g(\chi,\eta;k)
&=\chi\eta(-1)\tau(\eta)|k|^{1/2}H(1,k;\boldsymbol\chi,\boldsymbol\eta)
\end{align*}
with $\boldsymbol\chi=(\mathbf{1},\mathbf{1},\mathbf{1})$ and $\boldsymbol\eta=(\overline{\eta},\chi,\overline{\chi}),$ where $\mathbf{1}$ denotes the trivial character of $k^\times,$ $\tau(\eta)$ denotes the Gauss sum defined by $\eqref{eq:Gausssum}$ and $H$ denotes the hypergeometric sum defined by $\eqref{eq:hypergeometricsum}$.
\end{theorem}

The following inequality is an immediate consequence of Theorem \ref{thm:g-hyper} and Lemma \ref{lm:hypergeometric-weightranklisse} below.
\begin{corollary}
Let $\chi,\eta$ be two non-trivial characters of $k^\times.$ Then we have
\begin{align*}
|g(\chi,\eta;k)|\leqslant 3|k|.
\end{align*}
\end{corollary}

Our expression in Theorem \ref{thm:g-hyper} provides more possibilities for further studies on the (non-) correlations of $g(\chi,\eta)$ against other objects. For instance, as the application of a summation formula of Motohashi type, cubic moments of modular $L$-functions are related with the fourth moment of Dirichlet $L$-functions twisted by $g(\chi,\eta)$:
\begin{align*}
\sideset{}{^*}\sum_{\eta\bmod p}g(\chi,\eta) |L(\tfrac{1}{2},\eta)|^4.
\end{align*}
It is rather interesting to see if $g(\chi,\eta)$ correlates with $|L(\tfrac{1}{2},\eta)|^4$ as $\eta$ runs over all non-trivial multiplicative characters modulo $p$.

In a quite recent work \cite{WX23} on subconvex bounds for $\PGL_2$ $L$-functions over totally real number fields, we are led to study another double character sum
\begin{align*}
S(\chi,\eta;\rho):=\sum_{\alpha\in\bF_q} \rho(\alpha+\omega)\sum_{t\in\bF_q} \chi(t)\eta(\alpha^2-\omega^2 t)\overline{\eta}(1-t),
\end{align*}
where $\omega$ is a primitive element in $\bF_{q^2}$ such that $\omega^2\in\bF_q,$
$\rho$ is a non-trivial character of $\bF_{q^2}^{\times}$, and $\chi,\eta$ are both non-trivial characters of $\bF_q^{\times}$. By virtue of Gauss sums and Jacobi sums, we also transform
the inner sum over $t$ to a kind of hypergeometric sums, so that the upper bound for $S(\chi,\eta;\rho)$ follows from the non-correlation of such hypergeometric sum and the function $\alpha\mapsto \rho(\alpha+\omega),$ for which the theory of Deligne and Katz on Riemann Hypothesis for varieties over finite fields enters the picture.
It is worthy to mention that \cite{WX23} generalizes \cite{CI00} and \cite{PY20}, and in all works hypergeometric sums occupy important positions. As one may imagine, more features of hypergeometric sums are needed in \cite{WX23} due to the study of the (new) supercuspidal case therein. In view of the appearances of (hyper-) Kloosterman sums in Voronoi summation formulae, to which Motohashi's formula has the similar spirit, it is very acceptable that hypergeometric sums of Katz appear in the applications of Motohashi's formula and its generalizations; see \cite{CI00,PY20,WX23,Li11,Kw21} for instance.

\smallskip

\section{Gauss sums, Jacobi sums and hypergeometric sums}

The three objects mentioned in the title are intimately connected. For convenience, 
we write $\fG=\widehat{k^\times},$ the group of characters of $k^\times.$ For each $\chi\in\fG,$ it holds $\chi^{-1}=\overline{\chi}.$ We adopt the convention that $\chi(0)=0$ for each $\chi\in \fG.$ 

\subsection{Gauss sum}
For $\chi\in\fG$ and an additive character $\psi$ of $k,$ define the
Gauss sum
\begin{align}\label{eq:Gausssum}
\tau(\chi,\psi):= \sum_{\alpha\in\bF_q}\chi(\alpha)\psi(\alpha).
\end{align}
If $\psi$ is canonical, i.e., $\psi(x)=\exp(2\pi i \Tr_{k/\bF_p}(x)/p)$, we write $\tau(\chi,\psi)=\tau(\chi).$
It is an excise to show that
\begin{align*}
\tau(\chi,\psi)=-1
\end{align*}
if $\chi$ is trivial and $\psi$ is non-trivial.

\subsection{Jacobi sum}

For $\chi_1,\chi_2\in \fG,$ define the Jacobi sum
\begin{align*}
J(\chi_1,\chi_2):= \sum_{\alpha\in k}\chi_1(\alpha)\chi_2(1-\alpha).
\end{align*}

The following lemma associates Jacobi sums with Gauss sums, which is well-known, and should exist in literature for quite a long time; see Lidl and Niederreiter \cite[Theorem 5.21]{LN97} for instance.

\begin{lemma}\label{lm:Jacobi-Gauss}
Let $\chi_1,\chi_2\in \fG.$ Then $J(\chi_1,\chi_2)=|k|-2$ if $\chi_1,\chi_2$ are both trivial, and otherwise
\begin{align}\label{eq:Jacobi-Gauss}
J(\chi_1,\chi_2)= |k|^{-1}\tau(\chi_1,\psi)\tau(\chi_2,\psi)\overline{\tau(\chi_1\chi_2,\psi)}
\end{align}
for each non-trivial additive character $\psi$ of $k.$
\end{lemma}

\subsection{Hypergeometric sum}
We now introduce hypergeometric sums defined by Katz (see \cite[Chapter 8]{Ka90}). Let $m,n$ be two non-negative integers, and suppse $\boldsymbol{\chi}=(\chi_i)_{1\leqslant i\leqslant m}$ and $\boldsymbol{\eta}=(\eta_j)_{1\leqslant j\leqslant n}$ are two tuples of characters of $k^\times$, and $\psi$ is the canonical additive character of $k$. Katz introduced the following hypergeometric sum
\begin{align}\label{eq:hypergeometricsum}
H(t,k;\boldsymbol\chi,\boldsymbol\eta)
:=\frac{(-1)^{m+n-1}}{|k|^{(m+n-1)/2}}\mathop{\sum\sum}_{\substack{\bx\in(k^\times)^m,\by\in(k^\times)^n\\ \fn(\bx)=t\fn(\by)}}\boldsymbol\chi(\bx)
\overline{\boldsymbol\eta(\by)}\psi(\ft(\bx)-\ft(\by))
\end{align}
for $t\in k^\times,$ where, for $\bx=(x_1,x_2,\cdots,x_m)\in(k^\times)^m$, 
\begin{align*}
\boldsymbol\chi(\bx)=\prod_{1\leqslant i\leqslant m}\chi_i(x_i),
\end{align*}
\begin{align}\label{eq:T(x)N(x)}
\ft(\bx)=x_1+x_2+\cdots+x_m,\ \ \fn(\bx)=x_1x_2\cdots x_m,
\end{align}
and the notation with $\by$ can be defined in the same way. We say $\boldsymbol{\chi}$ and $\boldsymbol{\eta}$ are {\it disjoint} if $\chi_i\neq\eta_j$ for all $1\leqslant i\leqslant m$ and $1\leqslant j\leqslant n.$

\begin{lemma}\label{lm:hypergeometric-weightranklisse}
With the above notation,
if $\boldsymbol\chi$ and $\boldsymbol\eta$ are disjoint, then for any $\ell\neq p,$ there exists a geometrically irreducible $\ell$-adic middle-extension sheaf $\cH(\boldsymbol\chi,\boldsymbol\eta)$ on $\mathbf{A}_k^1$ with trace function given by $t\mapsto H(t,k;\boldsymbol\chi,\boldsymbol\eta),$ such that it is pointwise pure of weight zero and of rank $\max\{m,n\}.$
\end{lemma}

Lemma \ref{lm:hypergeometric-weightranklisse} was established by Katz \cite[Theorem 8.4.2]{Ka90}. One may refer to \cite{Ka90} for a very systematic study on geometric features of $H(t,q;\boldsymbol\chi,\boldsymbol\eta)$ and the underlying sheaf, which should be very crucial in further studies on the double character sum \eqref{eq:charactersum}. In particular, we may derive from Lemma \ref{lm:hypergeometric-weightranklisse} that
\begin{align*}
|H(t,k;\boldsymbol\chi,\boldsymbol\eta)|\leqslant \max\{m,n\}
\end{align*}
for all disjoint $\boldsymbol\chi$ and $\boldsymbol\eta$.

\smallskip

\section{Proof of Theorem \ref{thm:g-hyper}}

We first write
\begin{align*}
g(\chi,\eta;k)
&=\mathop{\sum\sum}_{u,v\in k^\times}\chi\Big(\frac{1+v^{-1}}{1+u^{-1}}\Big)\eta(uv-1).
\end{align*}
From the orthogonality of multiplicative characters, we have
\begin{align*}
g(\chi,\eta;k)
&=\sum_{x\in k}\eta(x-1)\mathop{\sum\sum}_{\substack{u,v\in k\\uvx=1}}\chi\Big(\frac{1+v}{1+u}\Big)\\
&=\frac{1}{|k|-1}\sum_{\rho\in \fG}\sum_{x\in k}\rho(x)\eta(x-1)\mathop{\sum\sum}_{u,v\in k}\rho(uv)\overline{\chi}(1-u)\chi(1-v)\\
&=\frac{\eta(-1)}{|k|-1}\sum_{\rho\in \fG} J(\rho,\eta)J(\rho,\overline{\chi})J(\rho,\chi).
\end{align*}

From Lemma \ref{lm:Jacobi-Gauss} it follows that
\begin{align*}
g(\chi,\eta;k)
&=\frac{\eta(-1)\tau(\eta)\tau(\overline{\chi})\tau(\chi)}{|k|^3(|k|-1)}\sum_{\rho\in\fG}\tau(\rho)^3\overline{\tau(\rho\eta)\tau(\rho\overline{\chi})\tau(\rho\chi)}.
\end{align*}
Opening the Gauss sums by definition, Theorem \ref{thm:g-hyper} then follows  immediately from orthogonality.

\smallskip

\subsection*{Acknowledgements} 
I would like to thank Yichen Qin and Han Wu for helpful conversations.
This work is supported in part by NSFC (No. 11971370, No. 12025106).

\smallskip

\bibliographystyle{plainnat}

\begin{thebibliography}{abcdefg}

\bibitem[CI00]{CI00} 
J. B. Conrey \& H. Iwaniec,
The cubic moment of central values of automorphic $L$-functions,
\emph{Ann. of Math. (2)} \textbf{151} (2000), 1175--1216. 


\bibitem[FKM21]{FKM21} 
\'{E}. Fouvry, E. Kowalski \& Ph. Michel, 
On the conductor of cohomological transforms,
\emph{Ann. Fac. Sci. Toulouse Math.} \textbf{30} (2021), no. 1, 203--254. 

\bibitem[Ka90]{Ka90} 
N. M. Katz,
Exponential Sums and Differential Equations, 
Ann. of Math. Stud., Vol. 124, Princeton University Press, Princeton, NJ, 1990.


\bibitem[Kw21]{Kw21} 
C.-H. Kwan,
Spectral moment formulae for $\GL(3)\times\GL(2)$ $L$-functions,
arXiv:2112.08568 [math.NT].


\bibitem[Li11]{Li11} 
X. Li,
Bounds for $GL(3)\times GL(2)$ $L$-functions and $GL(3)$ $L$-functions, \emph{Ann. of Math. (2)} \textbf{173} (2011), no. 1, 301--336.



\bibitem[LN97]{LN97} 
R. Lidl \& H. Niederreiter,
Finite Fields, 2nd Edition, Encyclopedia Math. Appl. \textbf{20}, Cambridge University Press, Cambridge, 1997.



\bibitem[PY20]{PY20} 
I. Petrow \& M.  Young,
The Weyl bound for Dirichlet $L$-functions of cube-free conductor,
\emph{Ann. of Math. (2)} \textbf{192} (2020), 437--486.

\bibitem[WX23]{WX23} 
H. Wu \& P. Xi,
A uniform Weyl bound for $L$-functions of Hilbert modular forms,
arXiv.

\end{thebibliography}

\end{document}